\documentclass[a4paper,english,12pt]{article}
\usepackage[logonly]{trace}
\usepackage{babel}
\usepackage{amsfonts, amsmath, amssymb}
\usepackage{graphicx}
\usepackage{pstricks}
\usepackage{psfig}
\usepackage{multirow}
\usepackage{fancyhdr}
\usepackage{vmargin, fancybox}

\newcommand{\mathsym}[1]{{}}

\usepackage{theorem} \theorembodyfont{\upshape}

\newtheorem{theorem}{Theorem}[section]

\newtheorem{definition}[theorem]{Definition}

\usepackage{graphicx}
\usepackage{pstricks}
\usepackage{psfig}

\makeatletter
\@addtoreset{equation}{section} \makeatother

\usepackage{amsmath}

\title{Lie group computation of finite difference schemes}

\begin{document}
\maketitle
\centerline{\scshape Emma Hoarau \footnotemark[1] and Claire David
\footnotemark[2]}
\medskip
\medskip
{\footnotesize \centerline{\footnotemark[1] ONERA, Computational
Fluid Dynamics and Aeroacoustics Department (DSNA)}
  \centerline{BP 72, 29 avenue de la Division Leclerc}
   \centerline{92322 Ch\text{\^ a}tillon Cedex, France}
\centerline{\footnotemark[2] Universit\'e Pierre et Marie
Curie-Paris 6}
  \centerline{Laboratoire de Mod\'elisation en M\'ecanique, UMR CNRS 7607}
   \centerline{Bo\^ite courrier $n^0 162$, 4 place Jussieu, 75252 Paris, cedex 05,
France}
}

\pagestyle{plain}

\medskip

\begin{quote}{\fontsize{8}{10}\selectfont
{\bfseries Abstract.} A \textsf{Mathematica} based program has been elaborated in order to determine the symmetry group of a finite difference equation.\\
The package provides functions which enable us to solve the determining equations of the related Lie group.
\par}
\end{quote}
\section{Introduction}
Various works have been carried out on the application of Lie group theory to numerical analysis. Most of them have been devoted to the building of numerical schemes, which preserve the symmetries of the original differential equations.\\
Olver \cite{Olvermf} and Kim \cite{Kim} constructed invariant finite difference equations using the concept of the moving frame. Discretization techniques, which preserve some symmetries of the original equations, were studied in \cite{Budd}, \cite{Dorodnitsyn}, \cite{Dorodnitsynal}, \cite{Bakirova} and \cite{Valiquette}. \\
\cite{Budd}, \cite{Dorodnitsyn}, \cite{Dorodnitsynal}, \cite{Bakirova} and \cite{Valiquette} used the discrete invariants of the Lie group of the original equation to build invariant finite difference schemes.\\
Yanenko \cite{Yanenko} and Shokin \cite{Shokin} have provided a Lie group analysis applied to finite difference equations, by means of a differential approximation and set down conditions under which the differential representation of a finite difference scheme preserves the symmetry group of the original differential equation. As the differential approximation is a differential equation, the Lie group theory can be fully applied to this notion.\\
The calculation of Lie groups of differential equations with pencil and paper is tedious and may induce errors. The size of related equations increases with the number of the symmetry variables, and the order of the differential equations. A large amount of packages have been created using software programs with symbolic manipulations, such as \textsf{Mathematica}, \textsf{MACSYMA}, \textsf{Maple}, \textsf{REDUCE}, \textsf{AXIOM}, \textsf{MuPAD}.
Schwarz \cite{Schwarz} wrote algorithms for \textsf{REDUCE} and \textsf{AXIOM} computer algebra systems, Vu and Carminati \cite{VuCarminati} worked on DESOLVE, a \textsf{Maple} program, Herod \cite{Herod} and Baumann \cite{Baumann} developed \textsf{Mathematica} programs.\\
We hereafter describe a new symbolic package, which implements the Lie group analysis methods for finite difference equations. The computations are based on the theory developed by Yanenko and Shokin. The program has been written for \textsf{Mathematica} and provides the symmetry group of a differential representation for a given finite difference scheme. The method is based on the \textsf{Mathematica} program of Cantwell in \cite{Cantwell}.\\
We presently aim at determining the symmetries lost by the discretization and building schemes which preserve those symmetries.\\
The method for the investigation of local point transformation groups is set out in section 2. The contents of the package is detailed in section 3. Implementation for classical numerical schemes is exposed in section 4.

\section{Lie group methods}

\subsection{Lie group of differential equations}

Consider a system of $l^{th}$-order differential equations:
\begin{equation}
\displaystyle \mathcal F^{\lambda}\big(x,u,u^{(k_1)},u^{(k_1,k_2)},\dots,u^{(k_1\dots k_l)}\big)=0,\ \ \lambda=1,\dots,q \label{eqn:ED}
\end{equation}
Denote by $u^{(k_1\dots k_p)}$ the vector, the components of which are partial derivatives of order $p$, namely, $u^{(k_1\dots k_p)}_j=\frac{\partial^p u_j}{\partial x_{k_1}\dots\partial x_{k_p}}$ $j=1,\dots,n$, $k_1,\dots, k_p \in \{1,\dots,m\}$.\\
Denote by $x=(x_1,\dots,x_m)$ the independent variables and $u=(u_1,\dots,u_n)$ the dependent variables.\\
The group of local point transformations can be written under the form:
\begin{equation}
\displaystyle G_r=\{x_i^{*}=\phi_i(x,u,a);\ u_j^{*}=\varphi_j(x,u,a),\ i=1,\dots,m;\ j=1,\dots,n\}
\label{eqn:ga}
\end{equation}

\noindent Expand the transformations by means of a Taylor series at the zero value of the parameter $a_\alpha$:
\begin{eqnarray}
x_i^{*}&=\displaystyle{x_i+a_\alpha\frac{\partial \phi_i}{\partial a_\alpha}\Big|_{a=0}+\mathcal{O}(a_\alpha^2),\ \alpha=1,\dots,r}\nonumber \\
u_j^{*}&=\displaystyle{u_j+a_\alpha\frac{\partial\varphi_j}{\partial a_\alpha}\Big|_{a=0}+\mathcal{O}(a_\alpha^2),\ \alpha=1,\dots,r}
\end{eqnarray}
The derivatives of $\phi_i$ and $\varphi_j$ with respect to the parameter $a_\alpha$ are smooth functions, called \emph{infinitesimals of the group} $G_r$. Denote by $\xi^\alpha_i$ and $\eta^\alpha_j$ the infinitesimals of $G_r$.

\noindent In order to find the Lie group transformations of the differential system, it is convenient to search the infinitesimal operators of $G_r$:
\begin{equation}
\displaystyle {\mathbf{L_\alpha}=\xi^\alpha_i(x,u)\frac{\partial}{\partial x_i}+\eta^\alpha_j(x,u) \frac{\partial }{\partial u_j},\ i=1,\dots,m;\ j=1,\dots,n;\ \alpha=1,\dots,r}
\label{eqn:opL}
\end{equation}
$\{\mathbf{L_\alpha}$, $\alpha=1,\dots,r\}$ represents the set of tangent vectors to the manifold $G_r$ at the neutral element $a=0$ and is a basis of the Lie-algebra of the infinitesimal operators of $G_r$.

\noindent The determination of the group transformations is reduced to the determination of the infinitesimal functions $\xi^\alpha_i$ and $\eta^\alpha_j$.

\noindent The knowledge of the $\mathbf{L_\alpha}$ enables us to determine the point transformations of the group $G_r$ by solving the equations:
\begin{eqnarray}
\displaystyle{\frac{\partial x_i^{*}}{\partial a_\alpha}}=\xi^\alpha_i(x^*,u^*),\ \displaystyle {\frac{\partial u_j^{*}}{\partial a_\alpha}}=\eta^\alpha_j(x^*,u^*),\ i=1,\dots,m;\ j=1,\dots,n;\ \alpha=1,\dots,r
\end{eqnarray}
in conjunction with the initial conditions:
\begin{eqnarray}
\displaystyle{x_i^{*}\big{|}_{a=0}}=x_i,\ \displaystyle{u_j^{*}\big{|}_{a=0}}=u_j
\end{eqnarray}

\vspace*{0.5cm}

\noindent In order to take into account the derivative terms involved in the differential equation, the Lie algebra vector field is prolonged:
\begin{eqnarray}
\displaystyle{\widetilde{\mathbf{L}}_\alpha^{(l)}=\mathbf{L_\alpha}+\sigma^{\alpha,(k_1)}_j
\frac{\partial}{\partial{u_j}^{(k_1)}}+\dots+\sigma^{\alpha,(k_1\dots k_l)}_j \frac{\partial}{\partial {u_j}^{(k_1\dots k_l)}}},\ \ \ \alpha=1,\dots,r;
\end{eqnarray}
$\sigma^{\alpha,(k_1)}_j$ and $\sigma^{\alpha,(k_1\dots k_o)}_j$ are given by:
\begin{eqnarray}
\sigma^{\alpha,(k_1)}_j &=& \displaystyle{\frac{\mathcal{D}\eta^\alpha_j}{\mathcal{D} x_{k_1}}-\sum_{i=1}^m \frac{\partial u_j}{\partial x_i} \frac{\mathcal{D}\xi^\alpha_i}{\mathcal{D}x_{k_1}}}\nonumber \\
\sigma^{\alpha,(k_1\dots k_o)}_j &=&\displaystyle{\frac{\mathcal{D}\sigma^{\alpha,(k_1\dots k_{o-1})}_j}{\mathcal{D} x_{k_o}^{\ \ \ \ }}-\sum_{i=1}^m \frac{\partial^{o} u_j}{\partial x_i \partial x_{k_1}\dots \partial
x_{k_{o-1}}} \frac{\mathcal{D}\xi^\alpha_i}{\mathcal{D}x_{k_o}},\ \ o=2,\dots,l}
\label{eqn:etas}
\end{eqnarray}
where:$\,\displaystyle\frac{\mathcal{D}}{\mathcal{D}x_k}=\frac{\partial}{\partial x_k}+\sum_{j=1}^{n}\frac{\partial u_j}{\partial x_k} \frac{\partial}{\partial u_j}$

\noindent Denote by $\widetilde{G}^{(l)}_r$ the Lie group of point transformations in the space\\ $\mathcal{E}\big(x,u,u^{(k_1)},u^{(k_1,k_2)},\dots,u^{(k_1\dots k_l)}\big)$ of the independent variables, the dependent variables and the derivatives of the dependent variables with respect to the independent ones.
\begin{definition}
Consider a subset $\Omega$ of the Euclidean space $\mathcal{E}\big(x,u,u^{(k_1)},u^{(k_1,k_2)},\dots,u^{(k_1\dots k_l)}\big)$
\begin{eqnarray}
\displaystyle \Omega=\{\big(x,u,u^{(k_1)},\dots,u^{(k_1\dots k_l)}\big); \mathcal F^{\lambda}\big(x,u,u^{(k_1)},\dots,u^{(k_1\dots k_l)}\big)=0,\ \lambda=1,\dots,q \}
\end{eqnarray}
$\Omega$ is an invariant subset of the group $\widetilde{G}^{(l)}_r$ if all the elements of $\widetilde{G}^{(l)}_r$ transform any point of $\Omega$ into a point of $\Omega$.
\end{definition}
\begin{theorem}
The system of $l^{th}$-order differential equations is invariant under the group $\widetilde G^{(l)}_r$ if and only if:
\begin{equation}
\displaystyle \widetilde{\mathbf{L}}_\alpha^{(l)}\mathcal F^\lambda{\Big{|}_{\Omega}} =0,\
\ \ \alpha=1,\dots,r;\ \lambda=1,\dots,q \label{eqn:deteq1}
\end{equation}
\label{th:invariance1}
\end{theorem}

\subsection{Lie group of differential approximations}
The finite difference scheme, which approximates the differential system (\ref{eqn:ED}), can be written as:
\begin{equation}
\displaystyle{\Lambda^{\lambda}(x,u,h,Tu)=0,\ \lambda=1,\dots,q}
\label{eqn:scheme}
\end{equation}
where $h=(h_1,h_2,\dots,h_m)$ denotes the space step vector, and $T=(T_1,T_2,\dots,T_m)$ the shift-operator along the axis of the independent variables, defined by:
\begin{equation}
T_i[u](x_1,x_2,\dots,x_{i-1},x_i,x_{i+1},\dots,x_m)=u(x_1,x_2,\dots,x_{i-1},x_i+h_i,x_{i+1},\dots,x_m).
\end{equation}

\begin{definition}
The differential equation:
\begin{eqnarray}
\displaystyle{\mathcal{P}^{\lambda}\big(x,u,u^{(k_1)},\dots,u^{(k_1\dots k_{l'})}\big)}&=& \displaystyle{\mathcal{F}^{\lambda}\big(x,u,u^{(k_1)},\dots,u^{(k_1\dots k_l)}\big)}\nonumber \\
& & +\displaystyle{\sum_{\beta=1}^s\sum_{i=1}^m (h_i)^{l_\beta} \mathcal{R}^{\lambda}_i(x,u,u^{(k_1)},\dots,u^{(k_1\dots k_{{l'}_{\lambda,i}})})},\nonumber \\
& & \displaystyle{\lambda=1,\dots,q};\ l'= max_{(\lambda,i)}
{l'}_{\lambda,i} \label{eqn:diffapprox}
\end{eqnarray}
is called the $l_s^{th}$-order differential approximation of the finite difference scheme (\ref{eqn:scheme}). In the specific case $s=1$, the above equation is called the first differential approximation.
\label{def:diffapprox}
\end{definition}

\noindent Denote by $G'_r$ a group of transformations in the space $\mathcal{E}(x,u,h)$:
\begin{equation}
\displaystyle G'_r=\{x_i^{*}=\phi_i(x,u,a);\
u_j^{*}=\varphi_j(x,u,a);h_i^{*}=\psi_i(x,u,h,a),\ i=1,\dots,m;\
j=1,\dots,n\}
\end{equation}
 by $\mathbf{L_\alpha}'$ the basis infinitesimal operator of $G'_r$:
\begin{equation}
\displaystyle{\mathbf{L_\alpha}'=
\mathbf{L_\alpha}+\zeta^\alpha_i(x,u,h)\frac{\partial}{\partial
h_i},\ \ \ \alpha=1,\dots,r}
\end{equation}
where
\begin{equation}
\displaystyle{\zeta^\alpha_i=\frac{\partial\psi_i}{\partial
a_\alpha}\Big{|}_{a=0},\ \ \alpha=1,\dots,r}
\end{equation}
and by $\widetilde{G}^{(l')}_r$ a group of transformation in the space $\mathcal{E}(x,u,h,u^{(k_1)},\dots,u^{(k_1\dots k_{l'})})$.\\
The ${l'}^{th}$-prolongation operator of $G'_r$, $\widetilde{\mathbf{L}}_\alpha^{(l')}$ can be written as:
\begin{equation}
\displaystyle
\widetilde{\mathbf{L}}_\alpha^{(l')}=\mathbf{L_\alpha}'+\sum_{j=1}^n\sum_{p=1}^{l'}\sigma_j^{\alpha,{(k_1\dots
k_p)}}\frac{\partial }{\partial u_j^{(k_1\dots k_p)}}
\label{eqn:prolongedop}
\end{equation}
\begin{theorem}
The differential approximation (\ref{eqn:diffapprox}) is invariant under the group $\widetilde G^{(l')}_r$ if and only if
\begin{equation}
\displaystyle{\widetilde{\mathbf{L}}_\alpha^{(l')}\mathcal{P}^{\lambda}\big((x,u,u^{(k_1)},\dots,u^{(k_1\dots
k_{l'})}\big)\Big|_{\mathcal{P}^{\lambda}=0}=0,\ \ \
\alpha=1,\dots,r;\ \lambda=1,\dots,q}
\end{equation}
or
\begin{equation}
\displaystyle{\Big[\widetilde{\mathbf{L}}_\alpha^{(l)}\mathcal{F}^{\lambda}+\widetilde{\mathbf{L}}_\alpha^{(l')}\Big(\sum_{\beta=1}^s\sum_{i=1}^m
(h_i)^{l_\beta}
\mathcal{R}^{\lambda}_i\Big)\Big]\Big|_{\mathcal{P}^{\lambda}=0}=0,\
\ \ \alpha=1,\dots,r;\ \lambda=1,\dots,q} \label{eqn:deteq2}
\end{equation}
\label{th:invariance2}
\end{theorem}

\noindent Equation (\ref{eqn:deteq2}) leads to a linear overdetermined system of partial differential equation, with respect to the infinitesimal functions, called the \textit{determining equations} of the Lie group of the differential approximation (\ref{eqn:diffapprox}). Our program determines the unknown infinitesimal functions.

\section{Computation Methods}

\noindent The program is restricted to partial differential approximations of any order involving the unknown scalar function $u$, the independent variables $x$ and $t$, the viscosity $\nu$ and the step size variables $h$ and $\tau$.

\subsection{Calculation of the differential approximation}
First, our program calculates the differential approximation from the knowledge of the considered finite difference equation and the approximation error.\\
The discrete approximations of the dependent variable, involved in the finite difference equation, are expanded at a given order by means of their Taylor series. The substitution of these Taylor series expansions into the finite difference scheme provides the $\Gamma$-form of the differential approximation, which contains derivatives with respect to $t$, $x$ and mixed derivatives with respect to $x$ and $t$. The $\Gamma$-form does not allow to have stability informations and does not yield the correct order of accuracy for all the numerical schemes. That's why we have found essential to determine directly the $\Pi$-form of the differential approximation, which is obtained by replacing the partial derivatives with respect to $t$ and mixed derivatives with respect to $x$ and $t$, involved in the $\Gamma$-form, by partial derivatives with respect to $x$, using the original differential equation.\\
The differential approximation is written under the form of an analytic function with respect to the independent variables, the dependent variable, the step size variables, the viscosity and the partial derivarives of the dependent variable. The dependent variable and its partial derivatives are considered as independent variables.

\subsection{Estimation of the determining equations}
Consider the infinitesimal functions $\xi^\alpha_i$,
$\eta^\alpha_j$, $\zeta^\alpha_i$ and $\chi^\alpha$. $\chi^\alpha$
is the infinitesimal related to the viscosity:
\begin{equation}
\displaystyle \chi^\alpha=\frac{\partial \pi}{\partial a_\alpha}\Big{|}_{a=0},\ \alpha=1,\dots,r
\end{equation}
where the viscosity transforms as follows $\nu^*=\pi(x,u,\nu,h,a)$.\\
The remaining infinitesimals of the prolongation operator of the considered symmetry group are generated according to the formulae (\ref{eqn:etas}).\\
The invariance condition of \textbf{Theorem \ref{th:invariance1}} provides a partial differential equation involving the unknown infinitesimal functions and products of the partial derivatives of the dependent variables.\\
Equation (\ref{eqn:deteq2}) is solved as an algebraic equation with respect to the partial derivatives of the dependent variables, handled as independent variables. Denote by $w$ the vector, the components of which are these variables. Since the whole equation holds for all the $w$ components, each coefficient in front of the products of the $w$ components has to be zero. This leads to a linear overdetermined system of partial differential equation, with respect to the infinitesimal functions, called the \textit{determining equations} of the Lie group of the differential system (\ref{eqn:diffapprox}). The overdetermined system is simplified by eliminating the redundancies. This step of the calculation requires the intervention of the user. The resolution of these equations yields explicitly the expression of $\xi^{\alpha}_i$, $\eta^{\alpha}_j$, $\zeta^\alpha_i$, $\chi^\alpha$, $\alpha=1,\dots,r,\ i=1,\dots,m,\ j=1,\dots,n$.

\subsection{Solving the determining equations}
The techniques used to solve the determining equations come from \cite{Cantwell}. The unknown infinitesimal functions are expanded by means of a power series expression with respect to the symmetry variables $x,\ t,\ u,\ \nu,\ h,\ \tau$. The polynomial expressions are substituted into the determining equations. Solving the determining equations amounts to finding the solutions of an algebraic equation.
Those latter techniques enable us to find in most cases the exact expression of the infinitesimals, when the power series is truncated, i.e. when the sought infinitesimals does not contain transcendental functions (like $\exp$, $\cos$, $\sin$, $\ln$, $\dots$).

\subsection{Determination of the symmetry group}

The last part of the program provides the infinitesimal function expression, the Lie algebra infinitesimal operators, and the corresponding Lie group transformations.

\section{Examples}
\noindent Consider the Burgers equation:
\begin{equation}
\displaystyle u_t+\frac{1}{2}(u^2)_x-\nu u_{xx}=0
\label{eqn:burgers}
\end{equation}
This equation admits the 6-parameter symmetry group:
\begin{itemize}
\item[\labelitemiii] $\displaystyle{\mathbf{L}_1=\frac{\partial}{\partial x}},\ \text{space translation}$

\item[\labelitemiii]  $\displaystyle{\mathbf{L}_2=\frac{\partial}{\partial t}},\ \text{time translation}$

\item[\labelitemiii]  $\displaystyle{\mathbf{L}_3=x\frac{\partial}{\partial x}+2 t\frac{\partial}{\partial t}-u\frac{\partial}{\partial u}},\ \text{dilatation}$
\item[\labelitemiii]  $\displaystyle{\mathbf{L}_4=x t\frac{\partial}{\partial x}+t^2\frac{\partial}{\partial t}+(-u t+x)\frac{\partial}{\partial u}},\ \text{projective transformation}$
\item[\labelitemiii]  $\displaystyle{\mathbf{L}_5=t\frac{\partial}{\partial
x}+\frac{\partial}{\partial u}},\ \text{Galilean transformation}$
\item[\labelitemiii]  $\displaystyle{\mathbf{L}_6=-t\frac{\partial}{\partial
t}+u\frac{\partial}{\partial u}+\nu\frac{\partial}{\partial
\nu}},\ \text{dilatation}$
\end{itemize}
Equation (\ref{eqn:burgers}) can be discretized by means of the finite difference schemes:
\begin{itemize}
\item[\labelitemiii]  \textbf{the FTCS (forward-time and centered-space) scheme}:
$$\displaystyle {\frac{u^{n+1}_i -u^n_i}{\tau}+\frac{\big(\frac{u^2}{2}\big)^n_{i+1} -\big(\frac{u^2}{2}\big)^n_{i-1}}{2 h}-\nu \frac{u^n_{i+1} -2 u^n_i+u^n_{i-1}}{h^2}=0}$$
\item[\labelitemiii]  \textbf{the Lax-Wendroff scheme}:
$$\displaystyle {\frac{u^{n+1}_i -u^n_i}{\tau}+\frac{\big(\frac{u^2}{2}\big)^n_{i+1} -\big(\frac{u^2}{2}\big)^n_{i-1}}{2 h}-\nu \frac{u^n_{i+1} -2 u^n_i+u^n_{i-1}}{h^2}+A^n_i=0}$$
\noindent where:
\begin{eqnarray}
\displaystyle A^n_i=&-&\frac{\tau}{2
h^2}\Big[u^n_{i+\frac{1}{2}}\Big(\big(\frac{u^2}{2}\big)^n_{i+1}
-\big(\frac{u^2}{2}\big)^n_{i}\Big)-u^n_{i-\frac{1}{2}}\Big(\big(\frac{u^2}{2}\big)^n_{i}-\big(\frac{u^2}{2}\big)^n_{i-1}\Big)\Big]\nonumber\\
&+&\frac{\nu
\tau}{2}\Big[\frac{\big(\frac{u^2}{2}\big)^n_{i+2}-2\big(\frac{u^2}{2}\big)^n_{i+1}
+2\big(\frac{u^2}{2}\big)^n_{i-1}-\big(\frac{u^2}{2}\big)^n_{i-2}}{h^3}\Big]\nonumber\\
&-&\frac{\nu^2\tau}{2}\Big[\frac{u^n_{i+2}-4 u^n_{i+1}+6 u^n_{i}-4
u^n_{i-1}+u^n_{i-2}}{h^4}\Big]\nonumber
\end{eqnarray}
\item[\labelitemiii]  \textbf{the Crank-Nicolson scheme}:
\begin{eqnarray*}
\displaystyle \frac{u^{n+1}_i-u^{n}_i}{\tau}&+&\frac{\big(\frac{u^2}{2}\big)^{n+1}_{i+1}-\big(\frac{u^2}{2}\big)^{n+1}_{i-1}+\big(\frac{u^2}{2}\big)^{n}_{i+1}-\big(\frac{u^2}{2}\big)^{n}_{i-1}}{4 h}\\
&-&\frac{\nu (u^{n+1}_{i+1}-2 u^{n+1}_i+u^{n+1}_{i-1}+u^{n}_{i+1} -2
u^{n}_i+u^{n}_{i-1})}{2 h^2}=0
\end{eqnarray*}
\end{itemize}
The first part of the program provides the differential representation of the schemes:
\begin{itemize}
\item[\labelitemiii]    \textbf{FTCS}
\begin{eqnarray*}
\displaystyle u_t+\frac{1}{2}(u^2)_x-\nu\ u_{xx}+\frac{\tau}{2}g_2+\frac{h^2}{12}(u^2)_{xxx}-\frac{\nu h^2}{12}u_{xxxx}=0
\end{eqnarray*}
\item[\labelitemiii]  \textbf{Lax-Wendroff}
\begin{eqnarray*}
\displaystyle u_t+\frac{1}{2}(u^2)_x-\nu\ u_{xx}+\frac{\tau^2}{6}g_3+\frac{h^2}{12}(u^2)_{xxx}-\frac{\nu h^2}{12}u_{xxxx}=0
\end{eqnarray*}
\item[\labelitemiii]  \textbf{Crank-Nicolson}
\begin{eqnarray*}
\displaystyle u_t+\frac{1}{2}(u^2)_x-\nu u_{xx}+\tau^2\Big(\frac{g_3}{6}+\frac{1}{4}(g_1^2+u g_2)_x-\frac{\nu}{4} (g_2)_{xx}\Big)+h^2\Big(\frac{1}{6}\big(\frac{u^2}{2}\big)_{xxx}-\frac{\nu}{12}u_{xxxx}\Big)=0
\end{eqnarray*}
where $g_1=-\big(\frac{u^2}{2}\big)_x+\nu u_{xx}$, $g_2=\big(-g_1 u\big)_x+\nu\big(g_1\big)_{xx}$, $g_3=\big(-g_2 u -g^2_1\big)_x+\nu \big(g_2\big)_{xx}$
\end{itemize}

\noindent The next steps of the computation of the symmetry group have been realized for the above schemes, but only the FTCS scheme is illustrated here.\\
The prolonged infinitesimal operator is calculated by means of formula (\ref{eqn:prolongedop}) and (\ref{eqn:etas}) with respect to the infinitesimal functions $\xi^\alpha_i$, $\eta^\alpha_j$, $\zeta^\alpha_i$ and $\chi^\alpha$.
The vector the components of which are the derivatives of the dependent variable, treated as independent variables, can be written as:
\begin{equation}
\displaystyle w=(wt,w2t,w3t,w4t,wx,wxt,wx2t,wx3t,w2x,w2xt,w2x2t,w3x,w3xt,w4x)
\end{equation}
The differential representation becomes:
\begin{equation}
\displaystyle wt+\frac{1}{2} u\ wx-\nu\ w2x+\frac{\tau}{2}\big(g_2^w\big)+h^2\big(\frac{1}{6}u\ w3x-\frac{1}{12}\nu\ w4x+\frac{1}{2} w2x\ wx\big)=0
\end{equation}
where $g^w_2=u^2\ w2x-2\nu\ u\ w3x+\nu^2\ w4x-4\nu\ w2x\ wx+ 2 u\ wx^2$.\\
The next step provides the determining equations, which are linear partial differential equations with respect to the unknown infinitesimal functions.
Some of them yield informations, which need to be entered by the user:
\begin{eqnarray}
\displaystyle & &\xi^\alpha_1=\xi^\alpha_1(x,t),\ \xi^\alpha_2=\xi^\alpha_2(t),\ \eta^\alpha_1=f_\eta(x,t)+u\ g_\eta(x,t),\nonumber\\[0.1cm]
\displaystyle & &\zeta^\alpha_1=\zeta^\alpha_1(x,t,u,h,\tau),\ \zeta^\alpha_2=\zeta^\alpha_2(x,t,u,h,\tau),\\[0.1cm]
\displaystyle & &\chi^\alpha=\chi^\alpha(x,t,u,h,\tau,\nu)\nonumber
\end{eqnarray}
Then the infinitesimal functions are expanded by means of multivariables polynomial expressions:
\begin{eqnarray}
\displaystyle & &\xi^\alpha_1=\sum^{\theta}_{i=0}\sum^{\theta-i}_{j=0} a_{ij}x^i y^j,\nonumber\\
\displaystyle & &\xi^\alpha_2=\sum^{\theta}_{i=0}b_i t^i,\nonumber\\[0.1cm]
\displaystyle & &\eta^\alpha_1=\sum^{\theta}_{i=0}\sum^{\theta-i}_{j=0} c_{ij}x^i y^j+u\sum^{\theta}_{i=0}\sum^{\theta-i}_{j=0} d_{ij}x^i y^j,\\[0.1cm]
\displaystyle & &\zeta^\alpha_1=\sum^{\theta}_{i=0}\Big(\sum^{\theta-i}_{j=0}\Big( \sum^{\theta-i-j}_{k=0}\Big(\sum^{\theta-i-j-k}_{l=0}\Big(\sum^{\theta-i-j-k-l}_{m=0} e_{ijklm}x^i y^j u^k h^l \tau^m\Big)\Big)\Big)\Big),\nonumber\\[0.1cm]
\displaystyle & &\zeta^\alpha_2=\sum^{\theta}_{i=0}\Big(\sum^{\theta-i}_{j=0}\Big( \sum^{\theta-i-j}_{k=0}\Big(\sum^{\theta-i-j-k}_{l=0}\Big(\sum^{\theta-i-j-k-l}_{m=0} f_{ijklm}x^i y^j u^k h^l \tau^m\Big)\Big)\Big)\Big),\nonumber\\[0.1cm]
\displaystyle & &\chi^\alpha=\sum^{\theta}_{i=0}\Big(\sum^{\theta-i}_{j=0}\Big( \sum^{\theta-i-j}_{k=0}\Big(\sum^{\theta-i-j-k}_{l=0}\Big(\sum^{\theta-i-j-k-l}_{m=0}\Big(\sum^{\theta-i-j-k-l-m}_{n=0} g_{ijklmn}x^i y^j u^k h^l \tau^m \nu^n\Big)\Big)\Big)\Big)\Big)\nonumber
\end{eqnarray}
The part performing the resolution of the algebraic determining equations gives the 4-parameter symmetry group represented by:
\begin{eqnarray}
\displaystyle \xi^\alpha_1&=&a_0+a_{10}\ x,\ \xi^\alpha_1=b_0+(2 a_{10}-g_1) t,\ \eta^\alpha_1=(g_1-a_{10})u,\ \alpha\leq4\\
\displaystyle \zeta^\alpha_1&=&a_{10} h,\ \zeta^\alpha_2=(2 a_{10}-g_1) h, \chi^\alpha=g_1 \nu\nonumber
\end{eqnarray}
The related 4-dimensional Lie algebra is generated by:
\begin{eqnarray}
\displaystyle L'_1&=&\frac{\partial}{\partial x},\ (a_0=1,b_0=0,a_{10}=0,g_1=0)\nonumber\\
\displaystyle L'_2&=&\frac{\partial}{\partial t},\ (a_0=0,b_0=1,a_{10}=0,g_1=0)\\
\displaystyle L'_3&=&x\frac{\partial}{\partial x}+2 t\frac{\partial}{\partial t}-u\frac{\partial}{\partial u}+h\frac{\partial}{\partial h}+2 \tau\frac{\partial}{\partial \tau},\ (a_0=0,b_0=0,a_{10}=1,g_1=0)\nonumber\\
\displaystyle L'_4&=&x\frac{\partial}{\partial x}+u\frac{\partial}{\partial u}+h\frac{\partial}{\partial h}+2\nu\frac{\partial}{\partial \nu},\ (a_0=0,b_0=0,a_{10}=0,g_1=1)\nonumber
\end{eqnarray}
The space and time translations preserve the differential representation of the FTCS scheme.\\
We recall that $(L_3,L_6)$ are the vectors of the basis of the Lie algebra $A_d$ of the Burgers equation dilatation group.\\
$(L'_3,L'_4)$ represent the vectors of the basis of the Lie algebra $A'_d$ of the differential approximation dilatation group. $A'_d$ can be expressed as a direct sum of the Lie algebra $A_{\alpha,\beta}$ spanned by the vectors \footnotesize{$\displaystyle (L_\alpha=x\frac{\partial}{\partial x}+2 t\frac{\partial}{\partial t}-u\frac{\partial}{\partial u},L_\beta=x\frac{\partial}{\partial x}+u\frac{\partial}{\partial u}+2\nu\frac{\partial}{\partial \nu})$} \normalsize and the Lie algebra $A_{h,\tau}$ spanned by the vectors \footnotesize{$\displaystyle(L'_\alpha=h\frac{\partial}{\partial h}+2 \tau\frac{\partial}{\partial \tau},L'_\beta=h\frac{\partial}{\partial h})$}\normalsize.\\
It is natural that if the independent and dependent variables are dilated then the step size variables undergo the same transformation.\\
$L_\alpha$ and $L_\beta$ are elements of the span of the set$\{L_3,L_6\}$. So $A'_d$ can be represented as the direct sum of $A_d$ and $A_{h,\tau}$. The FTCS scheme is invariant under the dilatation group, the Lie algebra of which can be written as $A_d \oplus A_{h,\tau}$.\\
The discretization by the FTCS scheme brings about the lost of the Galilean transformation and the projective transformation. The computation for the other schemes shows that all the finite difference schemes admit the same symmetry group.\\
The Lie group calculation program has been run under the \textsf{Mathematica 5.2} version for Windows on a PC with a Pentium 4, 2.6 GHz processor.
\begin{table}[!htp]
\begin{center}
\begin{tabular}{|c|c|c|}
    \hline  \textbf{Scheme} &   \textbf{differential approximation order}   &   \textbf{Time(s)}    \\
    \hline  FTCS            &   4                                           &   395.2               \\
    \hline  Lax-Wendroff    &   6                                           &   972.2               \\
    \hline  Crank-Nicolson  &   6                                           &   906.6               \\
    \hline
\end{tabular}
\end{center}
\caption{Characteristics}
\label{tabletime}
\end{table}
\noindent Table \ref{tabletime} shows the influence of the order of the differential approximation on the time of the calculation. The increase of differential approximation order yields a sharp rise of the time of the calculation.

\newpage

\end{document}